\documentclass[letterpaper, 10 pt, conference]{ieeeconf}  

\usepackage{amsmath}
\usepackage{graphicx}
\usepackage{cite}
\usepackage{balance}
\usepackage{amssymb,amsfonts}

\usepackage{amssymb,amsfonts}

\usepackage{lipsum}
\usepackage{mathtools}
\usepackage{cuted}

\usepackage{balance}

\DeclareMathOperator{\tr}{tr}

\newtheorem{theorem}{\it Theorem}
\newtheorem{lemma}{\it Lemma}
\newtheorem{definition}{\it Definition}

\newtheorem{corollary}{\it Corollary}

\IEEEoverridecommandlockouts                              

\title{\LARGE \bf
	The Spectral-Domain $\mathcal{W}_2$ Wasserstein Distance for Elliptical Processes and the Spectral-Domain Gelbrich Bound
}

\author{Song Fang$^{1}$ and Quanyan Zhu$^{1}$
\thanks{$^{1}$ Song Fang and Quanyan Zhu are with the Department of Electrical and Computer Engineering, New York University, New York, USA
        {\tt\small song.fang@nyu.edu; quanyan.zhu@nyu.edu}}%
}

\begin{document}

\maketitle
\thispagestyle{empty}
\pagestyle{empty}

\begin{abstract}
	
	In this short note, we introduce the spectral-domain $\mathcal{W}_2$ Wasserstein distance for elliptical stochastic processes in terms of their power spectra. We also introduce the spectral-domain Gelbrich bound for processes that are not necessarily elliptical.

\end{abstract}


\section{Introduction}

The Wasserstein distance (see, e.g., \cite{peyre2019computational, panaretos2020invitation} and the references therein) is an important metric from optimal transport theory (see, e.g., \cite{villani2003topics, villani2008optimal, santambrogio2015optimal} and the references therein). In this note, we consider the average $\mathcal{W}_2$ Wasserstein distance between stationary stochastic processes that are elliptical with the same density generator, which can be naturally characterized by a spectral-domain expression in terms of the power spectra of the processes (Theorem~\ref{t1}). On the other hand, when the stochastic processes are not necessarily elliptical, the average $\mathcal{W}_2$ Wasserstein distance is lower bounded by a spectral-domain generalization of the Gelbrich bound (Corollary~\ref{t2}).

\section{Preliminaries}

Throughout the note, we consider zero-mean real-valued continuous random variables and random vectors, as well as discrete-time stochastic processes. We represent random variables and random vectors using boldface letters, e.g., $\mathbf{x}$, while the probability density function of $\mathbf{x}$ is denoted as $p_\mathbf{x}$. We denote the sequence $\mathbf{x}_{0}, \ldots, \mathbf{x}_{k}$ by $\mathbf{x}_{0,\ldots,k}$ for simplicity, which, by a slight abuse of notation, is also identified with the random vector  $\left[ \mathbf{x}_0^{\mathrm{T}},\ldots,\mathbf{x}_{k}^{\mathrm{T}} \right]^{\mathrm{T}}$. 

We denote the covariance matrix of an $m$-dimensional random vector $\mathbf{x}$ by $\Sigma_{\mathbf{x}} = \mathbb{E}\left[  \mathbf{x} \mathbf{x}^{\mathrm{T}} \right]$. In the scalar case, the variance of $\mathbf{x}$ is denoted by $\sigma^2_{\mathbf{x}}$.
An $m$-dimensional stochastic process $\left\{ \mathbf{x}_{k}\right\}$ is said to be stationary if its auto-correlation $ R_{\mathbf{x}}\left( i,k\right)=\mathbb{E}\left[  \mathbf{x}_i \mathbf{x}_{i+k}^{\mathrm{T}} \right]$ depends only on $k$, and can thus be denoted as  $R_{\mathbf{x}}\left( k\right)$ for simplicity. The power spectrum of a stationary process $\left\{ \mathbf{x}_{k} \right\}$ is then defined as
\begin{flalign}
\Phi_{\mathbf{x}}\left( \omega\right)
=\sum_{k=-\infty}^{\infty} R_{\mathbf{x}}\left( k\right) \mathrm{e}^{-\mathrm{j}\omega k}. \nonumber
\end{flalign}
In the scalar case, the power spectrum of $\left\{ \mathbf{x}_{k} \right\}$ is denoted as $S_{\mathbf{x}}\left( \omega\right)$.
Note that throughout the note, all covariance matrices and power spectra are assumed to be positive definite.

The $\mathcal{W}_p$ Wasserstein distance (see, e.g., \cite{peyre2019computational, panaretos2020invitation}) is defined as follows.

\begin{definition} \label{Wdef}
	The $\mathcal{W}_p$ (for $p \geq 1$) Wasserstein distance between distribution $p_{\mathbf{x}}$ and distribution $p_{\mathbf{y}}$ is defined as
	\begin{flalign}
	\mathcal{W}_p \left( p_{\mathbf{x}} ; p_{\mathbf{y}} \right)
	= \left( \inf_{ \mathbf{x}, \mathbf{y}} \mathbb{E} \left[ \left\| \mathbf{x} - \mathbf{y} \right\|^p \right] \right)^{\frac{1}{p}}, \nonumber
	\end{flalign}
	where $\mathbf{x}$ and $\mathbf{y}$ denote $m$-dimensional random vectors with distributions $p_{\mathbf{x}}$ and $p_{\mathbf{y}}$, respectively.
\end{definition}

Particularly when $p = 2$, the $\mathcal{W}_2$ distance is given by 
\begin{flalign}
\mathcal{W}_2 \left( p_{\mathbf{x}} ; p_{\mathbf{y}} \right)
= \sqrt{ \inf_{ \mathbf{x}, \mathbf{y}} \mathbb{E} \left[ \left\| \mathbf{x} - \mathbf{y} \right\|^2 \right] }. \nonumber
\end{flalign}

%

The following lemma (see, e.g., \cite{ santambrogio2015optimal, peyre2019computational, panaretos2020invitation}) provides an explicit expression for the $\mathcal{W}_2$ distance between elliptical distributions with the same density generator. Note that Gaussian distributions are a special class of elliptical distributions (see, e.g., \cite{peyre2019computational}). Note also that hereinafter the random vectors are assumed to be zero-mean for simplicity.

\begin{lemma} \label{Gaussian}
	Consider $m$-dimensional elliptical random vectors $\mathbf{x}$ and $\mathbf{y}$ with the same density generator, while with covariance matrices $\Sigma_\mathbf{x}$ and $\Sigma_\mathbf{y}$, respectively.
	The $\mathcal{W}_2$ distance between distribution $p_\mathbf{x}$
	and distribution $p_\mathbf{y}$ is given by
	\begin{flalign}
	\mathcal{W}_2 \left( p_{\mathbf{x}} ; p_{\mathbf{y}} \right) 
	= \sqrt{ \tr \left[ \Sigma_{\mathbf{x}} + \Sigma_{\mathbf{y}} - 2 \left( \Sigma_{\mathbf{x}}^{\frac{1}{2}}  \Sigma_{\mathbf{y}}  \Sigma_{\mathbf{x}}^{\frac{1}{2}} \right)^{\frac{1}{2}} \right] }. \nonumber 
	\end{flalign}
\end{lemma}

\vspace{3mm}

Meanwhile, the Gelbrich bound (see, e.g., \cite{peyre2019computational}) is given as follows, which provides a generic lower bound for the $\mathcal{W}_2$ distance between distributions that are not necessarily elliptical.

\begin{lemma} \label{Gelbrich}
	Consider $m$-dimensional random vectors $\mathbf{x}$ and $\mathbf{y}$ with covariance matrices $\Sigma_\mathbf{x}$ and $\Sigma_\mathbf{y}$, respectively.
	The $\mathcal{W}_2$ distance between distribution $p_\mathbf{x}$
	and distribution $p_\mathbf{y}$ is lower bounded by
	\begin{flalign}
	\mathcal{W}_2 \left( p_{\mathbf{x}} ; p_{\mathbf{y}} \right) 
	\geq \sqrt{ \tr \left[ \Sigma_{\mathbf{x}} + \Sigma_{\mathbf{y}} - 2 \left( \Sigma_{\mathbf{x}}^{\frac{1}{2}}  \Sigma_{\mathbf{y}}  \Sigma_{\mathbf{x}}^{\frac{1}{2}} \right)^{\frac{1}{2}} \right] }. \nonumber
	\end{flalign}
\end{lemma}

\vspace{3mm}

\section{Spectral-Domain $\mathcal{W}_2$ Wasserstein Distance and Gelbrich Bound}

We first present the definition of average Wasserstein distance for stochastic processes.

\begin{definition}
	Consider $m$-dimensional stochastic processes $\left\{ \mathbf{x}_k \right\}$ and $\left\{ \mathbf{y}_k \right\}$.
	The average $\mathcal{W}_p$ distance between  $p_{\left\{ \mathbf{x}_k \right\}}$
	and $p_{\left\{ \mathbf{y}_k \right\}}$ is defined as
	\begin{flalign}
	&\mathcal{W}_p \left( p_{\left\{ \mathbf{x}_k \right\}} ; p_{\left\{ \mathbf{y}_k \right\}} \right) \nonumber \\ 
	&~~~~ = \left( \inf_{\left\{ \mathbf{x}_k \right\}, \left\{ \mathbf{y}_k \right\}} \limsup_{i \to \infty} \frac{ \mathbb{E} \left[ \left\| \mathbf{x}_{0,\ldots,i} - \mathbf{y}_{0,\ldots,i} \right\|^p \right]}{i+1} \right)^{\frac{1}{p}}.
	\end{flalign}
\end{definition}

\vspace{3mm}

In the case of $p = 2$, the average $\mathcal{W}_2$ distance is defined as \begin{flalign}
&\mathcal{W}_2 \left( p_{\left\{ \mathbf{x}_k \right\}} ; p_{\left\{ \mathbf{y}_k \right\}} \right) \nonumber \\ 
&~~~~ = \left( \inf_{\left\{ \mathbf{x}_k \right\}, \left\{ \mathbf{y}_k \right\}} \limsup_{i \to \infty} \frac{ \mathbb{E} \left[ \left\| \mathbf{x}_{0,\ldots,i} - \mathbf{y}_{0,\ldots,i} \right\|^2 \right]}{i+1} \right)^{\frac{1}{2}}.
\end{flalign}

We now proceed to present the main results of this note. The following Theorem~\ref{t1} provides a spectral-domain expression for the average $\mathcal{W}_2$ distance between elliptical processes.

\begin{theorem} \label{t1}
	Consider $m$-dimensional stationary stochastic processes $\left\{ \mathbf{x}_k \right\}$ and $\left\{ \mathbf{y}_k \right\}$ that are elliptical with the same density generator. Suppose that their distributions are given respectively by $p_{\left\{ \mathbf{x}_k \right\}}$ and $p_{\left\{ \mathbf{y}_k \right\}}$, while the power spectra are given respectively as $\Phi_\mathbf{x} \left( \omega \right)$ and $\Phi_\mathbf{y} \left( \omega \right)$.
	The average $\mathcal{W}_2$ distance between $p_{\left\{ \mathbf{x}_k \right\}}$ and $p_{\left\{ \mathbf{y}_k \right\}}$ is given by
	\begin{flalign} \label{W2m}
	&\mathcal{W}_2 \left( p_{\left\{ \mathbf{x}_k \right\}} ; p_{\left\{ \mathbf{y}_k \right\}} \right) \nonumber \\
	&~~~~ = \sqrt{ \frac{1}{2 \pi} \int_{0}^{2 \pi} \tr \left\{ W \left[ \Phi_\mathbf{x} \left( \omega \right), \Phi_\mathbf{y} \left( \omega \right) \right] \right\}  \mathrm{d} \omega },
	\end{flalign}
	where
	\begin{flalign} \label{W2m1}
	&W \left[ \Phi_\mathbf{x} \left( \omega \right), \Phi_\mathbf{y} \left( \omega \right) \right] \nonumber \\
	&~~~~ = \Phi_\mathbf{x} \left( \omega \right) + \Phi_\mathbf{y} \left( \omega \right) - 2 \left[ \Phi_\mathbf{x}^{\frac{1}{2}}  \left( \omega \right) \Phi_\mathbf{y} \left( \omega \right) \Phi_\mathbf{x}^{\frac{1}{2}} \left( \omega \right) \right]^{\frac{1}{2}}.
	\end{flalign}
\end{theorem}

\vspace{3mm}

\begin{proof}
	Note first that since $\left\{ \mathbf{x}_k \right\}$ and $\left\{ \mathbf{y}_k \right\}$ are stationary, we have
	\begin{flalign}
	&\mathcal{W}_2 \left( p_{\left\{ \mathbf{x}_k \right\}} ; p_{\left\{ \mathbf{y}_k \right\}} \right) \nonumber \\ 
	&~~~~ = \left( \inf_{\left\{ \mathbf{x}_k \right\}, \left\{ \mathbf{y}_k \right\}} \limsup_{i \to \infty} \frac{ \mathbb{E} \left[ \left\| \mathbf{x}_{0,\ldots,i} - \mathbf{y}_{0,\ldots,i} \right\|^2 \right]}{i+1} \right)^{\frac{1}{2}}
	\nonumber \\ 
	&~~~~ = \left( \inf_{\left\{ \mathbf{x}_k \right\}, \left\{ \mathbf{y}_k \right\}} \lim_{i \to \infty} \frac{ \mathbb{E} \left[ \left\| \mathbf{x}_{0,\ldots,i} - \mathbf{y}_{0,\ldots,i} \right\|^2 \right]}{i+1} \right)^{\frac{1}{2}}
	\nonumber \\ 
	&~~~~ = \left( \lim_{i \to \infty} \inf_{ \mathbf{x}_{0,\ldots,i}, \mathbf{y}_{0,\ldots,i} } \frac{  \mathbb{E} \left[ \left\| \mathbf{x}_{0,\ldots,i} - \mathbf{y}_{0,\ldots,i} \right\|^2 \right]}{i+1} \right)^{\frac{1}{2}}
	\nonumber \\ 
	&~~~~ = \left( \lim_{i \to \infty}  \frac{ \inf_{ \mathbf{x}_{0,\ldots,i}, \mathbf{y}_{0,\ldots,i} } \mathbb{E} \left[ \left\| \mathbf{x}_{0,\ldots,i} - \mathbf{y}_{0,\ldots,i} \right\|^2 \right]}{i+1} \right)^{\frac{1}{2}}.
	\nonumber
	\end{flalign}
	It then follows from Definition~\ref{Wdef} and Lemma~\ref{Gaussian} that 
	\begin{flalign}
	&\inf_{ \mathbf{x}_{0,\ldots,i}, \mathbf{y}_{0,\ldots,i} } \mathbb{E} \left[ \left\| \mathbf{x}_{0,\ldots,i} - \mathbf{y}_{0,\ldots,i} \right\|^2 \right] \nonumber \\
	& = \tr \left[ \Sigma_{\mathbf{x}_{0,\ldots,i}} + \Sigma_{\mathbf{y}_{0,\ldots,i}} - 2 \left( \Sigma_{\mathbf{x}_{0,\ldots,i}}^{\frac{1}{2}}  \Sigma_{\mathbf{y}_{0,\ldots,i}}  \Sigma_{\mathbf{x}_{0,\ldots,i}}^{\frac{1}{2}} \right)^{\frac{1}{2}} \right]
	\nonumber \\
	& = \tr \left( \Sigma_{\mathbf{x}_{0,\ldots,i}} \right) + \tr \left( \Sigma_{\mathbf{y}_{0,\ldots,i}} \right) \nonumber \\
	&~~~~~~~~ - 2 \tr \left[ \left( \Sigma_{\mathbf{x}_{0,\ldots,i}}^{\frac{1}{2}}  \Sigma_{\mathbf{y}_{0,\ldots,i}}  \Sigma_{\mathbf{x}_{0,\ldots,i}}^{\frac{1}{2}} \right)^{\frac{1}{2}} \right]
	, \nonumber
	\end{flalign}
	since $\mathbf{x}_{0,\ldots,i}$ and $\mathbf{y}_{0,\ldots,i}$ are elliptical with the same density generator. Meanwhile, it is known from \cite{gutierrez2008asymptotically, lindquist2015linear} that
	\begin{flalign}
	\lim_{i \to \infty} \frac{\tr \left( \Sigma_{\mathbf{x}_{0,\ldots,i}} \right)}{i+1}
	= \frac{1}{2 \pi} \int_{0}^{2 \pi} \tr \left[ \Phi_{\mathbf{x}} \left( \omega \right) \right] \mathrm{d} \omega, \nonumber
	\end{flalign}
	and
	\begin{flalign}
	\lim_{i \to \infty} \frac{\tr \left( \Sigma_{\mathbf{y}_{0,\ldots,i}} \right)}{i+1}
	= \frac{1}{2 \pi} \int_{0}^{2 \pi} \tr \left[ \Phi_{\mathbf{y}} \left( \omega \right) \right] \mathrm{d} \omega. \nonumber
	\end{flalign}
	It remains to prove that
	\begin{flalign}
	&\lim_{i \to \infty} \frac{\tr \left[ \left( \Sigma_{\mathbf{x}_{0,\ldots,i}}^{\frac{1}{2}}  \Sigma_{\mathbf{y}_{0,\ldots,i}}  \Sigma_{\mathbf{x}_{0,\ldots,i}}^{\frac{1}{2}} \right)^{\frac{1}{2}} \right]}{i+1} \nonumber \\
	&~~~~ = \frac{1}{2 \pi} \int_{0}^{2 \pi} \tr \left\{ \left[ \Phi^{\frac{1}{2}}_{\mathbf{x}} \left( \omega \right) \Phi_{\mathbf{y}} \left( \omega \right) \Phi^{\frac{1}{2}}_{\mathbf{x}} \left( \omega \right) \right]^{\frac{1}{2}} \right\} \mathrm{d} \omega. \nonumber
	\end{flalign}

	Note that $\Phi_{\mathbf{x}} \left( \omega \right)$ and $\Phi_{\mathbf{y}} \left( \omega \right)$ are positive definite. As such, $\Phi_\mathbf{x}^{\frac{1}{2}}  \left( \omega \right)$ is also positive definite (thus invertible), and
	\begin{flalign} 
	\Phi_\mathbf{x}^{\frac{1}{2}}  \left( \omega \right) \Phi_\mathbf{y} \left( \omega \right) \Phi_\mathbf{x}^{\frac{1}{2}} \left( \omega \right)
	=
	\Phi_\mathbf{x}^{- \frac{1}{2}}  \left( \omega \right) \left[ \Phi_\mathbf{x}  \left( \omega \right) \Phi_\mathbf{y} \left( \omega \right) \right] \Phi_\mathbf{x}^{ \frac{1}{2}}  \left( \omega \right), \nonumber
	\end{flalign}
	meaning that
	\begin{flalign} 
	\Phi_\mathbf{x}^{\frac{1}{2}}  \left( \omega \right) \Phi_\mathbf{y} \left( \omega \right) \Phi_\mathbf{x}^{\frac{1}{2}} \left( \omega \right) \nonumber
	\end{flalign}
	and 
	\begin{flalign} 
	\Phi_\mathbf{x} \left( \omega \right) \Phi_\mathbf{y} \left( \omega \right) \nonumber
	\end{flalign}
	are similar. Consequently, they share the same eigenvalues (which are all positive since $\Phi_\mathbf{x}^{\frac{1}{2}}  \left( \omega \right) \Phi_\mathbf{y} \left( \omega \right) \Phi_\mathbf{x}^{\frac{1}{2}} \left( \omega \right)$ is positive definite), while the square roots of the eigenvalues are also the same. Hence,
	\begin{flalign} 
	\tr \left[ \Phi_\mathbf{x}^{\frac{1}{2}}  \left( \omega \right) \Phi_\mathbf{y} \left( \omega \right) \Phi_\mathbf{x}^{\frac{1}{2}} \left( \omega \right) \right]^{\frac{1}{2}}
	= \tr \left[ \Phi_\mathbf{x} \left( \omega \right) \Phi_\mathbf{y} \left( \omega \right) \right]^{\frac{1}{2}}. \nonumber
	\end{flalign}
	Similarly, it can be proved that
	\begin{flalign}
	\tr \left[ \left( \Sigma_{\mathbf{x}_{0,\ldots,i}}^{\frac{1}{2}}  \Sigma_{\mathbf{y}_{0,\ldots,i}}  \Sigma_{\mathbf{x}_{0,\ldots,i}}^{\frac{1}{2}} \right)^{\frac{1}{2}} \right]
	= \tr \left[ \left( \Sigma_{\mathbf{x}_{0,\ldots,i}}  \Sigma_{\mathbf{y}_{0,\ldots,i}} \right)^{\frac{1}{2}} \right]. \nonumber
	\end{flalign}
	On the other hand, according to \cite{gutierrez2008asymptotically, lindquist2015linear}, we have
	\begin{flalign}
	&\lim_{i \to \infty} \frac{\tr \left[ \left( \Sigma_{\mathbf{x}_{0,\ldots,i}}  \Sigma_{\mathbf{y}_{0,\ldots,i}} \right)^{\frac{1}{2}} \right]}{i+1}
	\nonumber \\
	&~~~~ = \lim_{i \to \infty} \frac{ \sum_{j=0}^{\left(i + 1 \right)m - 1} \lambda_j \left[ \left( \Sigma_{\mathbf{x}_{0,\ldots,i}}  \Sigma_{\mathbf{y}_{0,\ldots,i}} \right)^{\frac{1}{2}} \right] }{i+1}
	\nonumber \\
	&~~~~ = \lim_{i \to \infty} \frac{ \sum_{j=0}^{\left(i + 1 \right)m - 1} \lambda^{\frac{1}{2}}_j  \left( \Sigma_{\mathbf{x}_{0,\ldots,i}}  \Sigma_{\mathbf{y}_{0,\ldots,i}} \right) }{i+1}
	\nonumber \\
	&~~~~= \frac{1}{2 \pi} \int_{0}^{2 \pi} \sum_{j=1}^{m} \lambda_{j}^{\frac{1}{2}}\left[ \Phi_{\mathbf{x}} \left( \omega \right) \Phi_{\mathbf{y}} \left( \omega \right) \right] \mathrm{d} \omega
	\nonumber \\
	&~~~~= \frac{1}{2 \pi} \int_{0}^{2 \pi} \sum_{j=1}^{m} \lambda_{j} \left\{ \left[ \Phi_{\mathbf{x}} \left( \omega \right) \Phi_{\mathbf{y}} \left( \omega \right) \right]^{\frac{1}{2}}  \right\} \mathrm{d} \omega
	\nonumber \\
	&~~~~= \frac{1}{2 \pi} \int_{0}^{2 \pi} \tr \left\{ \left[ \Phi_{\mathbf{x}} \left( \omega \right) \Phi_{\mathbf{y}} \left( \omega \right) \right]^{\frac{1}{2}}  \right\} \mathrm{d} \omega
	. \nonumber
	\end{flalign}
	Herein, we have used the fact that the eigenvalues of $\left( \Sigma_{\mathbf{x}_{0,\ldots,i}}  \Sigma_{\mathbf{y}_{0,\ldots,i}} \right)^{\frac{1}{2}}$ are given by the square roots of the eigenvalues of $\Sigma_{\mathbf{x}_{0,\ldots,i}}  \Sigma_{\mathbf{y}_{0,\ldots,i}}$, with the former being denoted as $
	\lambda_j \left[ \left( \Sigma_{\mathbf{x}_{0,\ldots,i}}  \Sigma_{\mathbf{y}_{0,\ldots,i}} \right)^{\frac{1}{2}} \right]$  while the latter as $
	\lambda_j \left( \Sigma_{\mathbf{x}_{0,\ldots,i}}  \Sigma_{\mathbf{y}_{0,\ldots,i}} \right)$.
	Similarly,
	the eigenvalues of $\left[ \Phi_{\mathbf{x}} \left( \omega \right) \Phi_{\mathbf{y}} \left( \omega \right) \right]^{\frac{1}{2}}$ are given by the square roots of the eigenvalues of $\Phi_{\mathbf{x}} \left( \omega \right) \Phi_{\mathbf{y}} \left( \omega \right)$, with the former being denoted as $\lambda_{j} \left\{ \left[ \Phi_{\mathbf{x}} \left( \omega \right) \Phi_{\mathbf{y}} \left( \omega \right) \right]^{\frac{1}{2}} \right\}$  while the latter as $
    \lambda_{j} \left[ \Phi_{\mathbf{x}} \left( \omega \right) \Phi_{\mathbf{y}} \left( \omega \right) \right]$.
	Consequently,
	\begin{flalign}
	&\lim_{i \to \infty} \frac{\tr \left[ \left( \Sigma_{\mathbf{x}_{0,\ldots,i}}^{\frac{1}{2}}  \Sigma_{\mathbf{y}_{0,\ldots,i}}  \Sigma_{\mathbf{x}_{0,\ldots,i}}^{\frac{1}{2}} \right)^{\frac{1}{2}} \right]}{i+1} \nonumber \\
	&~~~~ = \lim_{i \to \infty} \frac{\tr \left[ \left( \Sigma_{\mathbf{x}_{0,\ldots,i}}  \Sigma_{\mathbf{y}_{0,\ldots,i}} \right)^{\frac{1}{2}} \right]}{i+1}
	\nonumber \\
	&~~~~= \frac{1}{2 \pi} \int_{0}^{2 \pi} \tr \left\{ \left[ \Phi_{\mathbf{x}} \left( \omega \right) \Phi_{\mathbf{y}} \left( \omega \right) \right]^{\frac{1}{2}}  \right\} \mathrm{d} \omega \nonumber \\
	&~~~~= \frac{1}{2 \pi} \int_{0}^{2 \pi} \tr \left\{ \left[ \Phi^{\frac{1}{2}}_{\mathbf{x}} \left( \omega \right) \Phi_{\mathbf{y}} \left( \omega \right) \Phi^{\frac{1}{2}}_{\mathbf{x}} \left( \omega \right) \right]^{\frac{1}{2}} \right\}  \mathrm{d} \omega
	. \nonumber
	\end{flalign}
	This concludes the proof.
\end{proof}

Note that it is known from the proof of Theorem~\ref{t1} that
\begin{flalign} 
\tr \left[ \Phi_\mathbf{x}^{\frac{1}{2}}  \left( \omega \right) \Phi_\mathbf{y} \left( \omega \right) \Phi_\mathbf{x}^{\frac{1}{2}} \left( \omega \right) \right]^{\frac{1}{2}}
= \tr \left[ \Phi_\mathbf{x} \left( \omega \right) \Phi_\mathbf{y} \left( \omega \right) \right]^{\frac{1}{2}}.
\end{flalign}
Accordingly, $\tr \left\{ W \left[ \Phi_\mathbf{x} \left( \omega \right), \Phi_\mathbf{y} \left( \omega \right) \right] \right\}$ can equivalently be rewritten as
\begin{flalign} \label{W2m2}
&\tr \left\{ W \left[ \Phi_\mathbf{x} \left( \omega \right), \Phi_\mathbf{y} \left( \omega \right) \right] \right\} \nonumber \\
&~~~~ = \tr \left\{ \Phi_\mathbf{x} \left( \omega \right) + \Phi_\mathbf{y} \left( \omega \right) - 2 \left[ \Phi_\mathbf{x} \left( \omega \right) \Phi_\mathbf{y} \left( \omega \right) \right]^{\frac{1}{2}} \right\}.
\end{flalign}

Moreover, if it is further assumed that 
\begin{flalign}
	\Phi_\mathbf{x} \left( \omega \right) \Phi_\mathbf{y} \left( \omega \right) 
	= \Phi_\mathbf{y} \left( \omega \right) \Phi_\mathbf{x} \left( \omega \right),
\end{flalign}
then $W \left[ \Phi_\mathbf{x} \left( \omega \right), \Phi_\mathbf{y} \left( \omega \right) \right]$ reduces to
\begin{flalign} \label{W2m3}
W \left[ \Phi_\mathbf{x} \left( \omega \right), \Phi_\mathbf{y} \left( \omega \right) \right] = \Phi_\mathbf{x} \left( \omega \right) + \Phi_\mathbf{y} \left( \omega \right) - 2 \Phi_\mathbf{x}^{\frac{1}{2}} \left( \omega \right) \Phi_\mathbf{y}^{\frac{1}{2}} \left( \omega \right).
\end{flalign}
and accordingly, $\mathcal{W}_2 \left( p_{\left\{ \mathbf{x}_k \right\}} ; p_{\left\{ \mathbf{y}_k \right\}} \right)$ reduces to
\begin{flalign} 
&\mathcal{W}_2 \left( p_{\left\{ \mathbf{x}_k \right\}} ; p_{\left\{ \mathbf{y}_k \right\}} \right) \nonumber \\
& = \sqrt{ \frac{1}{2 \pi} \int_{0}^{2 \pi} \tr \left[ \Phi_\mathbf{x} \left( \omega \right) + \Phi_\mathbf{y} \left( \omega \right) - 2 \Phi_\mathbf{x}^{\frac{1}{2}} \left( \omega \right) \Phi_\mathbf{y}^{\frac{1}{2}} \left( \omega \right) \right]  \mathrm{d} \omega }
\nonumber \\
& = \sqrt{ \frac{1}{2 \pi} \int_{0}^{2 \pi} \left\| \Phi_\mathbf{x}^{\frac{1}{2}} \left( \omega \right) - \Phi_\mathbf{y}^{\frac{1}{2}} \left( \omega \right) \right\|_{\mathrm{F}}^2 \mathrm{d} \omega }.
\end{flalign}
This coincides with the Hellinger distance proposed in \cite{ferrante2008hellinger}. However, in general when $\Phi_\mathbf{x} \left( \omega \right)$ and $\Phi_\mathbf{y} \left( \omega \right)$ do not necessarily commute, it can be verified that
\begin{flalign} 
\mathcal{W}_2 \left( p_{\mathbf{x}} ; p_{\mathbf{y}} \right) 
\geq \sqrt{ \frac{1}{2 \pi} \int_{0}^{2 \pi} \left\| \Phi_\mathbf{x}^{\frac{1}{2}} \left( \omega \right) - \Phi_\mathbf{y}^{\frac{1}{2}} \left( \omega \right) \right\|_{\mathrm{F}}^2 \mathrm{d} \omega },
\end{flalign}
since
\begin{flalign} 
\tr \left[ \Phi_\mathbf{x}^{\frac{1}{2}}  \left( \omega \right) \Phi_\mathbf{y} \left( \omega \right) \Phi_\mathbf{x}^{\frac{1}{2}} \left( \omega \right) \right]^{\frac{1}{2}}
&\geq \tr \left[ \Phi_\mathbf{x}^{\frac{1}{4}} \left( \omega \right) \Phi_\mathbf{y}^{\frac{1}{2}} \left( \omega \right)
\Phi_\mathbf{x}^{\frac{1}{4}} \left( \omega \right) \right] \nonumber \\
& = \tr \left[ \Phi_\mathbf{x}^{\frac{1}{2}} \left( \omega \right) \Phi_\mathbf{y}^{\frac{1}{2}} \left( \omega \right) \right],
\end{flalign}
where the inequality is due to the Araki--Lieb--Thirring inequality \cite{araki1990inequality}; from this it also follows that
\begin{flalign} 
\tr \left\{ \left[ \Phi_\mathbf{x} \left( \omega \right) \Phi_\mathbf{y} \left( \omega \right)
 \right]^{\frac{1}{2}} \right\}
&= \tr \left[ \Phi_\mathbf{x}^{\frac{1}{2}}  \left( \omega \right) \Phi_\mathbf{y} \left( \omega \right) \Phi_\mathbf{x}^{\frac{1}{2}} \left( \omega \right) \right]^{\frac{1}{2}} \nonumber \\
& \geq  \tr \left[ \Phi_\mathbf{x}^{\frac{1}{2}} \left( \omega \right) \Phi_\mathbf{y}^{\frac{1}{2}} \left( \omega \right) \right],
\end{flalign}
which is a property that may be useful in other settings (herein the equality has been proved in the proof of Theorem~\ref{t1}).

Note also that in the scalar case (when $m = 1$), supposing that the power spectra of $\left\{ \mathbf{x}_k \right\}$ and $\left\{ \mathbf{y}_k \right\}$ are given respectively by $S_\mathbf{x} \left( \omega \right)$ and $S_\mathbf{y} \left( \omega \right)$,
	the average $\mathcal{W}_2$ distance between $p_{\left\{ \mathbf{x}_k \right\}}$
	and $p_{\left\{ \mathbf{y}_k \right\}}$ is given by
	\begin{flalign}
	&\mathcal{W}_2 \left( p_{\left\{ \mathbf{x}_k \right\}} ; p_{\left\{ \mathbf{y}_k \right\}} \right) \nonumber \\
	&= \sqrt{ \frac{1}{2 \pi} \int_{0}^{2 \pi} \left\{ S_\mathbf{x} \left( \omega \right) + S_\mathbf{y} \left( \omega \right) - 2 \left[ S_\mathbf{x} \left( \omega \right) S_\mathbf{y} \left( \omega \right) \right]^{\frac{1}{2}}  \right\} \mathrm{d} \omega } \nonumber \\
	&= \sqrt{ \frac{1}{2 \pi} \int_{0}^{2 \pi} \left[ S_\mathbf{x}^{\frac{1}{2}} \left( \omega \right) - S_\mathbf{y}^{\frac{1}{2}} \left( \omega \right) \right]^2 \mathrm{d} \omega}
	.
	\end{flalign}

The subsequent Corollary~\ref{t2} presents the spectral-domain Gelbrich bound for processes that are not necessarily elliptical, which follows directly from Theorem~\ref{t1}.

\begin{corollary} \label{t2}
	Consider $m$-dimensional stationary stochastic processes $\left\{ \mathbf{x}_k \right\}$ and $\left\{ \mathbf{y}_k \right\}$ that are not necessarily elliptical. Suppose that their distributions are given respectively by $p_{\left\{ \mathbf{x}_k \right\}}$ and $p_{\left\{ \mathbf{y}_k \right\}}$, while the power spectra are given respectively as $\Phi_\mathbf{x} \left( \omega \right)$ and $\Phi_\mathbf{y} \left( \omega \right)$.
	The average $\mathcal{W}_2$ distance between $p_{\left\{ \mathbf{x}_k \right\}}$ and $p_{\left\{ \mathbf{y}_k \right\}}$ is lower bounded by
	\begin{flalign} 
	&\mathcal{W}_2 \left( p_{\left\{ \mathbf{x}_k \right\}} ; p_{\left\{ \mathbf{y}_k \right\}} \right) \nonumber \\
	&~~~~ \geq \sqrt{ \frac{1}{2 \pi} \int_{0}^{2 \pi} \tr \left\{ W \left[ \Phi_\mathbf{x} \left( \omega \right), \Phi_\mathbf{y} \left( \omega \right) \right] \right\}  \mathrm{d} \omega },
	\end{flalign}
	where
	\begin{flalign} 
	&W \left[ \Phi_\mathbf{x} \left( \omega \right), \Phi_\mathbf{y} \left( \omega \right) \right] \nonumber \\
	&~~~~ = \Phi_\mathbf{x} \left( \omega \right) + \Phi_\mathbf{y} \left( \omega \right) - 2 \left[ \Phi_\mathbf{x}^{\frac{1}{2}}  \left( \omega \right) \Phi_\mathbf{y} \left( \omega \right) \Phi_\mathbf{x}^{\frac{1}{2}} \left( \omega \right) \right]^{\frac{1}{2}}.
	\end{flalign}
\end{corollary}

\vspace{3mm}

In the scalar case, the spectral-domain Gelbrich bound is given by
\begin{flalign}
&\mathcal{W}_2 \left( p_{\left\{ \mathbf{x}_k \right\}} ; p_{\left\{ \mathbf{y}_k \right\}} \right) \nonumber \\
&\geq \sqrt{ \frac{1}{2 \pi} \int_{0}^{2 \pi} \left\{ S_\mathbf{x} \left( \omega \right) + S_\mathbf{y} \left( \omega \right) - 2 \left[ S_\mathbf{x} \left( \omega \right) S_\mathbf{y} \left( \omega \right) \right]^{\frac{1}{2}}  \right\} \mathrm{d} \omega } \nonumber \\
&= \sqrt{ \frac{1}{2 \pi} \int_{0}^{2 \pi} \left[ S_\mathbf{x}^{\frac{1}{2}} \left( \omega \right) - S_\mathbf{y}^{\frac{1}{2}} \left( \omega \right) \right]^2 \mathrm{d} \omega}
.
\end{flalign}

\section{Conclusion}

	In this note, we have introduced the spectral-domain $\mathcal{W}_2$ Wasserstein distance for elliptical stochastic processes. We have also introduced the spectral-domain Gelbrich bound for processes that are not necessarily elliptical. It might be interesting to examine the implications of the results in future (cf. \cite{W2properties} for instance).

\balance

\bibliographystyle{IEEEtran}
\bibliography{references}

\begin{thebibliography}{10}
\providecommand{\url}[1]{#1}
\csname url@samestyle\endcsname
\providecommand{\newblock}{\relax}
\providecommand{\bibinfo}[2]{#2}
\providecommand{\BIBentrySTDinterwordspacing}{\spaceskip=0pt\relax}
\providecommand{\BIBentryALTinterwordstretchfactor}{4}
\providecommand{\BIBentryALTinterwordspacing}{\spaceskip=\fontdimen2\font plus
\BIBentryALTinterwordstretchfactor\fontdimen3\font minus
  \fontdimen4\font\relax}
\providecommand{\BIBforeignlanguage}[2]{{%
\expandafter\ifx\csname l@#1\endcsname\relax
\typeout{** WARNING: IEEEtran.bst: No hyphenation pattern has been}%
\typeout{** loaded for the language `#1'. Using the pattern for}%
\typeout{** the default language instead.}%
\else
\language=\csname l@#1\endcsname
\fi
#2}}
\providecommand{\BIBdecl}{\relax}
\BIBdecl

\bibitem{peyre2019computational}
G.~Peyr{\'e} and M.~Cuturi, ``Computational optimal transport: With
  applications to data science,'' \emph{Foundations and Trends{\textregistered}
  in Machine Learning}, vol.~11, no. 5-6, pp. 355--607, 2019.

\bibitem{panaretos2020invitation}
V.~M. Panaretos and Y.~Zemel, \emph{An Invitation to Statistics in Wasserstein
  Space}.\hskip 1em plus 0.5em minus 0.4em\relax Springer, 2020.

\bibitem{villani2003topics}
C.~Villani, \emph{Topics in Optimal Transportation}.\hskip 1em plus 0.5em minus
  0.4em\relax American Mathematical Society, 2003.

\bibitem{villani2008optimal}
------, \emph{Optimal Transport: Old and New}.\hskip 1em plus 0.5em minus
  0.4em\relax Springer, 2008.

\bibitem{santambrogio2015optimal}
F.~Santambrogio, \emph{Optimal Transport for Applied Mathematicians}.\hskip 1em
  plus 0.5em minus 0.4em\relax Birk{\"a}user, Springer, 2015.

\bibitem{gutierrez2008asymptotically}
J.~Guti{\'e}rrez-Guti{\'e}rrez and P.~M. Crespo, ``Asymptotically equivalent
  sequences of matrices and {H}ermitian block {T}oeplitz matrices with
  continuous symbols: Applications to {MIMO} systems,'' \emph{IEEE Transactions
  on Information Theory}, vol.~54, no.~12, pp. 5671--5680, 2008.

\bibitem{lindquist2015linear}
A.~Lindquist and G.~Picci, \emph{Linear Stochastic Systems: A Geometric
  Approach to Modeling, Estimation and Identification}.\hskip 1em plus 0.5em
  minus 0.4em\relax Springer, 2015.

\bibitem{ferrante2008hellinger}
A.~Ferrante, M.~Pavon, and F.~Ramponi, ``Hellinger versus {K}ullback--{L}eibler
  multivariable spectrum approximation,'' \emph{IEEE Transactions on Automatic
  Control}, vol.~53, no.~4, pp. 954--967, 2008.

\bibitem{araki1990inequality}
H.~Araki, ``On an inequality of {L}ieb and {T}hirring,'' \emph{Letters in
  Mathematical Physics}, vol.~19, no.~2, pp. 167--170, 1990.

\bibitem{W2properties}
S.~Fang and Q.~Zhu, ``Independent elliptical distributions minimize their
  $\mathcal{W}_2$ {W}asserstein distance from independent elliptical
  distributions with the same density generator,'' \emph{arXiv preprint}, 2020.

\end{thebibliography}

\end{document}